\documentclass{amsart}

\newtheorem{theorem}{Theorem}[section]
\newtheorem{lemma}[theorem]{Lemma}
\newtheorem{prop}[theorem]{Proposition}

\theoremstyle{definition}

\theoremstyle{remark}

\numberwithin{equation}{section}

\begin{document}

\title{Infinite Product Representations for Multiple Gamma Function}

\author{Michitomo Nishizawa}

\address{Graduate School of Mathematical Science,
            the University of Tokyo,
         3-8-1 Komaba, Meguro, Tokyo, 153-8914 Japan.}
\email{mnishi@ms.u-tokyo.ac.jp}

\thanks{This research is partially supported by Grant-in-Aid for Scientific Research
No.1407387 from Japan Society for the Promotion of Science.}





\begin{abstract}
Two kinds of infinite product representations 
for Vign\'eras multiple gamma function are presented. 
As an application of these formulas,
a multiplication formula for the function is derived. 
\end{abstract}

\maketitle

\section*{Introduction}

In a series of papers \cite{bar1, bar2, bar3, bar4},  
Barnes introduced multiple gamma functions  
associated with a certain generalization of the Hurwitz zeta function.
In relevant with a special case of Barnes' function,
Vign\'eras \cite{vig} introduced her multiple gamma functions $G_r(z)$  
$(r \in \mathbb{Z}_{\geq0})$ as a sequence of meromorphic
functions uniquely determined by the following relations:
\begin{align}
 & \mbox{(i)}\ G_0(z) = z, \quad  
   \mbox{(ii)}\ G_r(1) = 1, \quad
   \mbox{(iii)}\ G_r(z+1) = G_{r-1}(z)G_r(z) \nonumber\\[-8pt]
 & \label{gBM1} \\[-8pt]
 & \mbox{(iv)}\ \frac{d^{r+1}}{dz^{r+1}} \log G_r(z+1) \geq 0 \quad
   \mbox{for}\ z \geq 0. \nonumber
\end{align}
This formulation can be considered as a generalization of 
the Bohr-Morellup theorem.
For example, $G_1(z)$ is the celebrated Euler gamma function
$\Gamma(z)$ ({\it cf}. Artin \cite{art}, Whittaker-Watson \cite{whi}.).
$G_2(z)$ is $G$-function introduced in Barnes \cite{bar1}.\par
Vign\'eras multiple gamma function has various applications to number theory,
geometry and mathematical physics.  
It is known that Vign\'eras' function appears in factors in the
determinants of the Laplacians on the some compact 
Riemann surfaces.
Sarnak \cite{sar} applied  $G_2(z)$ to the representation of the determinants
of the Laplacian on spinor fields on a Riemann surface.
Many researchers \cite{cho, CS, kum, var, vor}  
computed the determinants in the case where the Laplacians are on the $n$-sphere. 
$G_2(z)$ is used to represent factors of T\"oplitz determinants.
For example, we can refer to B\"ottcher-Silbermann \cite{BS}. 
Tracy \cite{tra} and Basor-Tracy \cite{BT} applied the fact
to a representation of coefficients in an asymptotic expansion for
$\tau$-function of the Ising model. 
As mentioned in \cite{ada1}, it is expected to apply Vign\'eras'
function to random matrix theory .
$G_2(z)$ also appears in an asymptotic behavior of the mean value of a certain
$L$-function \cite{CF}.
Vign\'eras' functions are closely related with Kurokawa's 
multiple sine function \cite{kur1, kur2, kur3, KK}. The gamma factor
of the Selberg zeta function is represented as products of these
functions.\par
As its significance is gradually recognized,  
studies on Vign\'eras' function have become of a much interest 
to researchers in the theory of special functions. 
Ueno and the author \cite{UN} derived an asymptotic expansion
({\it ``higher Stirling formula''})
and  an Weierstrass' infinite product representation for $G_r(z)$. 
The author \cite{nis2} gave an another proof of Barnes'
generalization of the H\"older theorem \cite{bar5}.
It is proved that Vign\'eras functions
satisfy no algebraic differential equation. 
We should note results by Ferreira and Lopez \cite{FL} 
and by Adamtik \cite{ada1, ada2}.\par
In this paper, we present two types of infinite product representations 
of Vign\'eras' multiple gamma function, 
which can be considered as a generalization of 
the Gauss and of the Euler product formula of Euler's gamma function
\begin{align}
 \Gamma(z+1) & = \lim_{N\to \infty} \frac{N!}{(z+1)(z+2)\cdots(z+N)}
     (N+1)^z \label{gamma1}\\[4pt]
    & = \prod_{n=1}^{\infty} \left[
        \left(1+\frac{z}{n}\right)^{-1} \left(1+\frac1n\right)^z
      \right]\label{gamma2}
\end{align} 
({\it cf.}  Artin \cite{art}, Whittaker-Watson \cite{whi}).
Our main theorem is stated as follows:
If $z$ is not negative integer,  the multiple gamma function $G_r(z)$
is represented as
\begin{align}
  G_r (z+1) & = \lim_{N\to \infty} \left[
    \prod_{n=1}^{N} \frac{G_{r-1}(n)}{G_{r-1}(z+n)}
    \prod_{k=0}^{r-1} G_{k} (N+1)^{\binom{z}{r-k}}
  \right] \label{intGauss}\\[4pt]
   & = \prod_{n=1}^{\infty} \left[
      \frac{G_{r-1}(n)}{G_{r-1}(z+n)} 
      \prod_{k=0}^{r-1}\left(
         \frac{G_k(n+1)}{G_k(n)}
      \right)^{\binom{z}{r-k}}
       \right]. \label{intEuler}
\end{align} 
In the case when $r=1$, these formulas coincide with (\ref{gamma1}) and
(\ref{gamma2}).  We can find the representation for $G_2(z)$ in Jackson \cite{jac}.
It should be noted that infinite product formula of these types for 
a $q$-analogue of the multiple gamma function 
were already obtained in \cite{nis1}.
However, in contrast to simplicity in $q$-case, 
some delicate techniques are necessary
to deal with infinite products of Vign\'eras' function.
We verify (\ref{intGauss}) and (\ref{intEuler})
in section 1.
The point is to apply an asymptotic expansion in \cite{UN}
to estimations for products of Vign\'eras' functions.
\par
In section 2, as an application of infinite product representations, 
we derive a multiplication formula for Vign\'eras' multiple gamma function,
which can be regarded as a generalization of the well known formula 
\begin{equation}
   \prod_{m=0}^{p}\Gamma\left( \frac{z+m}{p}\right)
    = \frac{(2\pi)^{\frac{p-1}2}}{p^{z-\frac12}}\Gamma(z)
\label{Emulp}\end{equation}
for Euler's gamma function ({\it cf.}  Artin \cite{art}, Whittaker-Watson \cite{whi}).
It is described as follows:
\[
   \prod_{q_1,\ q_2, \cdots q_r=0}^{p-1} 
   G_r \left(
     \frac{z+q_1+\cdots + q_r}{p}
   \right)
   = \frac{e^{\phi_r(z)}}{p^{\psi_r(z)}} G(z)
\]
It might be seem that formula of this type can be guessed easily from (\ref{gBM1}).
However, it is not easy to determine explicit forms of $\phi_r(z)$ and
$\psi_r(z)$. The reason why we can do it is usefulness of 
our representations (\ref{intGauss}). \par
On preparing the first draft of this paper, 
the author was informed about Kuribaysahi's result \cite{kuri}.
He proved a multiplication formula 
\[
  \prod_{q_1, \cdots, q_r =0}^{p-1}
   \Gamma_r \left(
     \frac{z+q_1+\cdots +q_r}p
   \right)
    = p^{Q_r(z)} \Gamma_r(z).
\]
for a function $\Gamma_r(z)$ defined as
\[
   \Gamma_r(z):= \exp  \left[\left.
     \frac{\partial}{\partial s}
     \zeta_r(s,z)
    \right|_{s=0}\right],
\]
where $\zeta_r(s,z)$ is a generalization of the Hurwitz zeta function
defined as the analytical continuation of the series
\[
  \zeta_r(s,z) := \sum_{n_1,\cdots , n_r=0}^{\infty}
     (z+n_1+\cdots +n_r)^{-s}, \qquad
  \Re s > r.
\]
At the end of section 2, 
discussions about relations between his result and ours are added. 
We can find the following relation:
\begin{equation}
 Q_r(z) = (-1)^r \psi_r(z),
\label{IntQPsi}\end{equation}
although they look different at first sight.
It may be worth noting that the relation (\ref{IntQPsi}) 
seems to be applicable to Kurokawa's multiple sine function 
\cite{kur1, kur2, kur3, KK}.
(\ref{IntQPsi}) can be verified without use of the zeta
function. We give an elementary proof in appendix.\par
For simplicity, we call Vign\'eras multiple gamma function 
only ``multiple gamma function'' in the following sections. \\[12pt]
\hspace{-4pt}{\bf Notations}:\
In this paper, we use notation $B_r(z)$
for the Bernoulli polynomial defined by the generating function
\[
   \sum_{r=0}^{\infty} B_r(z) t^r
    = \frac{te^{t}}{1-e^t},
\] 
and $B_r$ for the Bernoulli number defined as $B_r:=B_r(0)$.
We introduce the Stirling number ${}_r S_j$ 
of the 1 st kind by 
\[
  t (t-1) \cdots (t-r+1) = \sum_{j=0}^r {}_r S_j t^j.   
\] 
The notation $\zeta(s)$ is used to refer to the Riemann zeta function
defined as the series  
$\zeta(s):=\sum_{n=1}^{\infty} n^{-s}$
and its analytical continuation. $\zeta'(s)$ is the first derivative 
of $\zeta(s)$ defined by $\zeta'(s) := \frac{d}{ds}\zeta(s)$. 

\section{Infinite Product Representation}
As mentioned in introduction, our main theorem is described as
follows:
\begin{theorem}
If $z$ is not negative integer and is included in any finite region 
of complex plane,  the multiple gamma function $G_r(z)$
is represented as
{\allowdisplaybreaks
\begin{align}
  G_r (z+1) & = \lim_{N\to \infty} \left[
    \prod_{n=1}^{N} \frac{G_{r-1}(n)}{G_{r-1}(z+n)}
    \prod_{k=0}^{r-1} G_{k} (N+1)^{\binom{z}{r-k}}
  \right] \label{Gauss}\\[4pt]
   & = \prod_{n=1}^{\infty} \left[
      \frac{G_{r-1}(n)}{G_{r-1}(z+n)} 
      \prod_{k=0}^{r-1}\left(
         \frac{G_k(n+1)}{G_k(n)}
      \right)^{\binom{z}{r-k}}
       \right].\label{Euler}
\end{align} 
}
\label{main}\end{theorem}
\noindent
\begin{proof} 
From the Gauss product representation (\ref{Gauss}), 
the Euler product representation (\ref{Euler})
follows immediately. 
So, we give a proof of (\ref{Gauss}) in this section.
We apply an asymptotic expansion for $G_r(z)$,
which was firstly appeared in \cite{UN}.

\begin{theorem}[Ueno-Nishizawa]
Let us put $0<\delta<\pi$, then, 
as $|z|\to\infty$ in the sector $\{z \in {\bold C}||\arg z|<\pi - 
\delta\}$,
\begin{align}
& \log G_r(z+1) = \left\{
    \binom{z+1}{r} +\sum_{j=0}^{r-1} \frac{B_{j+1}}{j+1} G_{r,j}(z)
  \right\} \log(z+1) - \nonumber\\[-8pt]
& \label{hS}\\[-8pt]
& \qquad - \sum_{j=0}^{r-1} G_{r,j}(z) \frac{(z+1)^{j+1}}{(j+1)^2}
    - \sum_{j=0}^{r-1} G_{r,j}(z) \zeta'(-j) +O(z^{-1}).\nonumber
\end{align}
where a polynomial $G_{r,j}(z)$ is defined by the generating function
\[
   \binom{z-u}{r-1} =: \sum_{j=0}^{r-1} G_{r,j}(z)u^j
   \quad (r=0 \cdots r-1), \qquad
   G_{r,j}(z) = 0, \quad (j \geq r).
\] 
\label{UN}\end{theorem}
\noindent
In our proof, the following lemma is useful:
\begin{lemma}
For arbitrary $x$, $y \in \mathbb{C}$,
\begin{align*}
 & {\rm (i)}\quad \sum_{k=0}^r \binom{x}{r-k} \binom{y}{k}
        = \binom{x+y}{r}, \qquad
  {\rm (ii)}\quad \sum_{k=0}^{r} \binom{x}{r-k} G_{k,j} (y) =
   G_{r,j}(x+y). 
\end{align*}
\label{binom}\end{lemma}

\noindent
Noting this lemma and that
\[
   \sum_{j=0}^{r-1} G_{r,j}(z+N-1) \left\{
     \frac{(z+N)^{j+1}}{j+1} - \frac{N^{j+1}}{(j+1)^2}
   \right\}
    = \int_N^{z+N} \frac{dv}{v} \int_0^v \binom{z+N-1-u}{r-1} du,
\]
we rewrite the logarithms of terms in brackets of (\ref{Gauss})
and have the following asymptotic behavior as $N\to\infty$:
{\allowdisplaybreaks
\begin{align*}
& \log  \left[
    \prod_{n=1}^{N} \frac{G_{r-1}(n)}{G_{r-1}(z+n)}
    \prod_{k=0}^{r-1} G_{k} (N+1)^{\binom{z}{r-k}}
  \right] =\\
&  =   \log G_r(z+1) + 
     \sum_{k=0}^{r} \binom{z}{r-k} 
     \left\{
       \binom{N}{k} + \sum_{j=0}^{r-1} \frac{B_{j+1}}{j+1} G_{k,j}(N-1)
     \right\} \log N - \\ 
& - \sum_{j=0}^{r-1} \sum_{k=0}^{r} \binom{z}{r-k} G_{k,j}(N-1)
    \frac{N^{j+1}}{j+1} 
    - \sum_{j=0}^{r-1} \sum_{k=0}^{r} \binom{z}{r-k} G_{k,j}(N-1)
    \zeta'(-j) - \\
& - \left\{
     \binom{z+N}{r} + \sum_{j=0}^{r-1} \frac{B_{j+1}}{j+1} G_{r,j}(z+N-1)
   \right\} \log (z+N) + \\
& +   \sum_{j=0}^{r-1} G_{r,j}(z+N-1) \frac{(z+N)^{j+1}}{j+1}
   +   \sum_{j=0}^{r-1} G_{r,j}(z+N-1)\zeta'(-j) + O(N^{-1}) = \\
& = \log G_r(z+1) + \int_0^z \frac{du}{u+N} \left[
   \binom{z+N}{r}+ \sum_{j=1}^r \frac{B_{j+1}}{j+1} G_{r,j-1}(z+N-1)
   - \right.\\
& \qquad - \left. \int_{-1}^{z+N-1} \binom{v}{r-1} dv
   + \int_u^z \binom{z-1-v}{r-1}dv 
  \right] +O(N^{-1}).
\end{align*}}
\noindent
As $N\to\infty$, this integral vanishes 
because of the following lemma:
\begin{lemma}[Ueno-Nishizawa]
For arbitrary $z \in \mathbb{C}$, we have
\[
  \binom{z}{r} + \sum_{j=0}^{r-1} \frac{B_{j+1}}{j+1} G_{r,j}(z-1)
   = \int_{-1}^{z-1} \binom{t}{r-1} dt
     + \sum_{j=0}^{r-1} \frac{B_{j+1}}{j+1} G_{r,j}(-1).
\]
\label{UN36}\end{lemma}
\noindent
This was already shown in \cite{UN}.
Therefore, we have proved theorem \ref{main}. 
\end{proof}

\section{Multiplication formula}
As an application of Gauss' product representation, we demonstrate the 
multiplication formula of the multiple gamma function. 

\begin{theorem}
\begin{equation}
   \prod_{q_1,\ q_2, \cdots q_r=0}^{p-1} 
   G_r \left(
     \frac{z+q_1+\cdots + q_r}{p}
   \right)
   = \frac{e^{\phi_r(z)}}{p^{\psi_r(z)}} G(z)
\label{Vmulp}\end{equation}
where
\begin{align*}
& \phi_r(z) = \sum_{j=0}^{r-1} \left[
    \sum_{q_1,\cdots, q_r=0}^{p-1}
    G_{r,j} \left(
      \frac{z+q_1+\cdots+q_r}p -2 
    \right)
    -G_{r,j} \left(z-1\right)\right] \zeta'(-j)
   \\
& \psi_r(z) = \binom{z}{r}+\sum_{j=0}^{r-1} \frac{B_{j+1}}{j+1} G_{r,j}(z-1).
\end{align*}
\label{mult}\end{theorem}
\begin{proof}
From the infinite product representation (\ref{Gauss}), 
it follows that
\begin{align*}
& \prod_{q_1,\cdots, q_r=0}^{p-1}
   G_r \left(
    \frac{z+q_1+\cdots +q_r}p
   \right) =  \\
& \qquad = \lim_{N\to\infty}\left[
     \frac{\prod_{n=1}^{p(N-1)-1} G_{r-1}(n)}
       {\prod_{n=1}^{p(N-1)-1} G_{r-1}(z+n-1)}
     \times \prod_{k=0}^{r-1} G_r(p(N-1))^{\binom{z-1}{r-k}}
   \right] \times \\
& \qquad \times \lim_{N\to\infty}\left[
   \prod_{k=0}^r 
    \frac{G_r (N)^{\sum_{q_1,\cdots, q_r}\binom{(z+q_1+\cdots+q_r)/p-1}{r-k}}}
      {G_r(p(N-1))^{\binom{z-1}{r-k}}} 
    \times 
    \frac{p^{\sum_{m=0}^{pN-1} \psi_{r-1}(z+m) }}
     {e^{\sum_{m=0}^{pN-1} \phi_{r-1}(z+m)}}
  \right]
\end{align*}
We substitute the asymptotic expansion (\ref{hS}) to the logarithm of
terms in the second bracket. 
{\allowdisplaybreaks
\begin{align*}
 & \log \left[
   \prod_{k=0}^r 
    \frac{G_r (N)^{\sum_{q_1,\cdots, q_r}\binom{(z+q_1+\cdots+q_r)/p-1}{r-k}}}
      {G_r(p(N-1))^{\binom{z-1}{r-k}}} 
    \times 
    \frac{p^{\sum_{m=0}^{pN-1} \psi_{r-1}(z+m) }}
     {e^{\sum_{m=0}^{pN-1} \phi_{r-1}(z+m)}}
  \right]= \\ 
 & = \left\{
    \sum_{q_1,\cdots, q_r}\binom{(z+q_1+\cdots+q_r)/p-1}{r}
    + \sum_{j=0}^{r-1} \frac{B_{j+1}}{j+1} G_{r,j} \left(
       \frac{z+q_1+\cdots+q_r}p+N-2
     \right)
    \right\} \log N - \\
 & - \left\{ \binom{z+p(N-1)-1}{r}
    + \sum_{j=0}^{r-1} \frac{B_{j+1}}{j+1} G_{r,j} \left(
       (z+p(N-1)-2)
     \right) 
    \right\} \log \left(
      N-1-\frac1p
     \right) - \\
 &  - \sum_{j=0}^{r+1}\left[
     \sum_{q_1,\cdots,q_r}  G_{r,j} \left(
         \frac{z+q_1+\cdots+q_r}p+N-2
       \right)  -
       G_{r,j}(z+p(N-1)-1)
    \right]\frac{N^{j+1}}{(j+1)^2} - \\
 &  + \left\{
    \binom{z+p(N-1-1)}{r} + \sum_{j=0}^{r-1} \frac{B_{j+1}}{j+1} G_{r,j}
    (z+p(N-1)-1) - \sum_{m=0}^{p(N-1)-1} \psi_r(z+m)
   \right\} \log p - \\
 &  -  \sum_{j=0}^{r-1}\left[
     \sum_{q_1,\cdots,q_r=0}^{p-1}
       G_{r,j} \left(
         \frac{z+q_1+\cdots+q_r}p+N-2
       \right)  -
       G_{r,j}(z+p(N-1)-1)
    \right]\zeta'(-j) - \\
 & \qquad - \sum_{m=0}^{pN-1} \phi_r (z+m) + o(1).
\end{align*}
}
\hspace{-4pt}We show that its divergent terms vanish.
First, we compute terms including $\log p$.
\begin{prop}
If we define $\psi_0(z)=0$ and
 \[
   \psi_r (z) := \binom{z}{r} + \sum_{j=0}^{r-1} \frac{B_{j+1}}{j+1} G_{r,j}(z-1),
 \]
then $\psi_r(z)$ satisfies $\psi_0(z)=z$ and 
\[
  \binom{p(N-1)-1}{r} + \sum_{j=0}^{r-1} \frac{B_{j+1}}{j+1} G_{r,j}(z+p(N-1)-2)
   - \sum_{m=0}^{p(N-1)-1} \psi_{r-1}(z+m) = \psi_r(z).
\] 
$\psi_r(z)$ does not depend on $N$ and is uniquely determined  as the polynomial satisfying the
above recurrence relation. 
\end{prop}
\noindent
\begin{proof} 
This proposition immediately follows from the relation 
\[
  \sum_{l=0}^{L-1} \binom{x+k}{k} = \binom{z+L}{k+1} - \binom{z}{k+1}
\]
for $L\in \mathbb{Z}_{\geq 0}$.
\end{proof}
\noindent
Next, we simplify terms including $\zeta'(-j)$
and give a explicit form of $\phi_r(z)$.
\begin{prop}
If we define 
\[
  \phi_{r,j} (z) := \sum_{q_1, \cdots, q_r=0}^{p-1}
   G_{r,j} \left(
    \frac{z+q_1+\cdots+q_r}p - 2
    \right)
   - G_{r,j} (z-1),
\]
then 
$
  \phi_r(z) = \sum_{j=0}^{r-1} \phi_{r,j}(z) \zeta'(-j)
$
is uniquely determined as a polynomial satisfying
the recurrence relation $\phi_0(z)=0$ and  
\begin{align*}
 & \sum_{j=0}^{r-1}\left[
     \sum_{q_1,\cdots,q_r=0}^{p-1}
       G_{r,j} \left(
         \frac{z+q_1+\cdots+q_r}p+N-2
       \right)  -
       G_{r,j}(z+p(N-1)-1)
    \right]\zeta'(-j) - \\
 & \qquad - \sum_{m=0}^{pN-1} \phi_{r-1} (z+m) = \phi_r(z)
\end{align*}
\end{prop}
\noindent
\begin{proof}
It is sufficient to prove
\begin{align*}
& \phi_{r,j}(z) = \sum_{q_1,\cdots, q_r} 
    G_{r,j} \left(
      \frac{z+q_1+\cdots+q_r}p +N-2
    \right) - G_{r,j}(z+p(N-1)-1) - \\
&- \sum_{m=0}^{p(N-1)-1}\left[
     \sum_{q_1,\cdots, q_r} 
    G_{r-1,j} \left(
      \frac{z+m+q_1+\cdots+q_{r-1}}p +N-2
    \right) - G_{r-1,j}(z+m+p(N-1)-1)
  \right] .
\end{align*}
However, we can conclude it from the identity
\[
  \sum_{m=0}^{L} G_{r,j}(z+m) = G_{r+1,j}(z+L) - G_{r+1,j}(z),
  \quad (L\in \mathbb{Z}_{\geq 0}),
\]
and
\begin{align*}
&  \sum_{m=0}^{p(N-1)-1} \sum_{q_1,\cdots, q_{r-1}=0}^{p-1}
   G_{r,j} \left(
     \frac{z+m+q_1+\cdots +q_{r-1}}p-2
   \right)=\\
& = \sum_{q_1, \cdots, q_{r-1}, q=0}^{p-1}
 \left[
  G_r \left(
   \frac{z+q_1+\cdots +q_r}p + N -2
  \right) -G_r \left(
 \frac{z+q_1+\cdots +q_r}p -2
  \right)
  \right],
\end{align*}
The uniqueness of $\phi_r(z)$ follows from its polynomiality.
\end{proof}
In order to finish our proof, we verify that
the rest of terms  vanish as $N\to\infty$.
By lemma \ref{binom}, we can see that
{\allowdisplaybreaks
\begin{align*}
& \sum_{k=0}^r \left\{
    \sum_{q_1, \cdots, q_r=0}^{p-1}
    \binom{(z+q_1+\cdots +q_r)/p - 1}{r-k}
  \right\} \times \\
& \times  \left[
    \binom{N}{k} + \sum_{j=0}^{k-1} \frac{B_{j+1}}{j+1} G_{k,j} (N-1)
     - \sum_{j=0}^r G_{k,j} (N-1) \frac{N^2}{(j+1)^2}
  \right] - \\
& - \sum_{k=0}^r \binom{z-1}{r-k} \left\{
    \binom{N}{k} + \sum_{j=0}^{r-1} \frac{B_{j+1}}{j+1} G_{k,j} (N-1)
    -\sum_{j=0}^{k-1} G_{k,j} (N-1) \frac{N^{j+1}} {(j+1)^2}
  \right\} \\
& = \sum_{q_1, \cdots, q_r=0}^{p-1} \left[
    \binom{(z+q_1+\cdots +q_r)/p - 1}{r} + 
    \sum_{j=0}^{r-1} \frac{B_{j+1}}{j+1} G_{r,j} \left(
      \frac{z+q_1+\cdots +q_r}p+N-2
    \right)- \right.\\
& \left.   - \sum_{j=0}^{r-1} G_{r,j} \left(
      \frac{z+q_1+\cdots +q_r}p+N-2
    \right) \frac{N^2}{(j+1)^2}     
   \right] - \\
& - \left\{
   \binom{z+N-1}{r} 
  + \sum_{j=0}^{r-1} \frac{B_{j+1}}{j+1} G_{r,j}(z+N-2) 
  - \sum_{j=0}^{r-1} G_{r,j}(z+N-2)\frac{N^{j+1}}{(j+1)^2}
 \right\}.
\end{align*}
}
\hspace{-4pt}From the same argument as proof of theorem \ref{main},
it follows that the above terms tend to zero as $N\to\infty$. 
Therefore, we have proved theorem \ref{mult}.
\end{proof}
Our result is closely related with Kuribayashi \cite{kuri}.
In order to explain his result, we introduce some functions.
$\zeta_r(s,z)$ is defined as a special case of Barnes' zeta function 
\cite{bar4, var}, which is introduced as the series 
\[
  \zeta_r (s,z) := \sum_{n_1,\cdots, n_r=0}^{\infty}
     (z+n_1+\cdots + n_r)^{-s} 
\]
for $\Re s > r$. This function can be continued analytically 
to a meromorphic function whose poles are placed at
$s=1, \cdots ,r$.
We call the analytic continuation also $\zeta_r(s,z)$.
The gamma function $\Gamma_r(z)$ associated with $\zeta_r(s,z)$
is introduced as
\[
 \Gamma_r(z):=\exp\left[\left.
    \frac{\partial}{\partial s} \zeta_r (s,z)
   \right|_{s=0} \right].
\]
Kuribayashi exhibit the following multiplication formula:

\begin{theorem}[Kuribayashi]
$\Gamma_r(z)$ satisfies the following multiplication formula:
\[
  \prod_{q_1, \cdots, q_r=0}^{p-1} \Gamma_r \left(
   \frac{z+q_1+\cdots + q_r}p
  \right) 
  = p^{Q_r(z)} \Gamma_r (z),
\]
where
\[
  Q_r(z) = \frac{(-1)^r}{(r-1)!}\sum_{r=1}^{r} \frac{{}_r S_l}l \left\{
      z^l - (-1)^l B_l
    \right\}.
\]
\label{kuri}\end{theorem}
\noindent
As a consequence of facts in Vardi \cite{var},
a relation between $G_r (z)$ and $\Gamma_r(z)$ is expressed as follows:
\[
  G_r(z) = R_r(z) \Gamma_r(z)^{(-1)^{r-1}}
\]
where 
\begin{align*}
&  R_r(z) := \exp\left[
     \sum_{j=0}^{r-1} G_{r,j} (z-1) \zeta'(-j)
     \right]
\end{align*} 
Thus, we have
\begin{equation}
  Q_r(z) = (-1)^r \psi_r(z) 
    =  (-1)^r \left[
     \binom{z}{r} + \sum_{j=0}^{r-1}
      \frac{B_{j+1}}{j+1} G_{r,j} (z-1)
   \right].
\label{Qpsi}\end{equation}
Our expression is useful in some cases of studies on related functions.
For example, noting that $G_{r,0}(z) = \binom{z}{r-1}$,
we can check that the relation follows 
\begin{equation}
 (-1)^r Q_r(r-z) = Q_r (z).
\label{hanten}\end{equation}
from the definition of $\psi_r(z)$ and (\ref{Qpsi}). 
It plays an important role in the multiplication
formula 
\[
  \prod_{q_1,\cdots , q_r=0}^{p-1}
    S_r \left(
     \frac{z+q_1+\cdots +q_r}p
    \right)
    = S_r(z).
\]
for Kurokawa's multiple sine function
\cite{kur1, kur2, kur3, KK} introduced as
\[
   S_r(z) := \Gamma_r(r-z)\Gamma_r(z)^{(-1)^{r+1}}.
\]
In Kuribayashi's original proof,  (\ref{hanten})
is verified through a rather complicated argument,
He applied a relation
between $\zeta_r(-m,z)$ ($m\in \mathbb{Z}_{\geq 0}$)
and the Bernoulli polynomials $B_l(z)$.
However, once (\ref{Qpsi}) is obtained,
we can check (\ref{hanten}) immediately.

\section{Appendix : An Elementary Proof for (\ref{Qpsi})}
Without facts of zeta functions, we can prove (\ref{Qpsi}) 
directly as follows:
First, we rewrite Kuribayashi's $Q_r(z)$ as 
\begin{align}
& (-1)^r Q_r(z) = \frac1{(r-1)!} \sum_{l=0}^{r-1}
    {}_{r-1} S_l \left\{
     \frac{(-1)^{l+1} B_{l+1}}{l+1}
      - \frac{(z-1)^{l+1}}{l+1}
    \right\}. 
\label{modQ}\end{align}
The second term can be written as follows:
\begin{align*}
& \frac1{(r-1)!} \sum_{l=0}^{r-1} 
    {}_{r-1} S_l \frac{(z-1)^{l+1}}{l+1} 
  = \int_0^z \binom{t-1}{r-1} dt
    - \int_0^1 \binom{t-1}{r-1} dt.
\end{align*}
From Lemma \ref{UN36} and 
\[
  G_{r,j} (0) = \frac{(-1)^j}{(r-1)!} {}_{r-1} S_j,
\]
it follows that
\begin{align*}
& \int_0^z \binom{t-1}{r-1} dt
    - \int_0^1 \binom{t-1}{r-1} dt = \\
& \quad = \binom{z}{r} 
    + \sum_{j=0}^{r-1} \frac{B_{j+1}}{j+1} G_{r,j}(z-1)
   - \frac1{(r-1)!} \sum_{j=0}^{r-1} \frac{B_{j+1}}{j+1}(-1)^j {}_{r-1}S_{j}. 
\end{align*}
Therefore, we obtain (\ref{Qpsi}) by substituting this to (\ref{modQ}). 
\section*{acknowledgements}
The author expresses his deep gratitude to Professor Michio Jimbo 
and Professor Kimio Ueno  for their encouragement. 
He is also grateful to  Professor Nobushige Kurokawa 
for his stimulating lecture at Waseda University.
He also thanks to Masanori Kuribayashi for informing his result \cite{kuri}
and to Yuji Hara and Yashushi Kajihara for their interests to this work. 

\bibliographystyle{amsplain}

\end{document}